\newcommand{\R}{I\!\!R}      % insiemi numerici
\newcommand{\Z}{\mathbb Z}   %%%%%%%%%%%%%%%%%%

\newcommand{\N}{I\!\!N}
\newcommand{\Spl}{\mathrm{Sp}}
\newcommand{\spl}{\mathrm{sp}}

\newcommand{\iMaslov}{\mathrm i_{\textrm{Maslov}}}   % indici vari
\newcommand{\Scal}{\mathcal S}
\newcommand{\Kcal}{\mathcal K}
\newcommand{\Bcal}{\mathcal B}
\newcommand{\Qcal}{\mathcal Q}
\newcommand{\Fcal}{\mathcal F}
\newcommand{\Fcalsa}{\mathcal F^{\scriptscriptstyle{\textrm{sa}}}}
\newcommand{\ess}{{\!\scriptscriptstyle{\textrm{ess}}}}
\newcommand{\spfl}{\mathrm{sf}}
\newcommand{\Hcal}{\mathcal H}

\newcommand{\Ker}{\mathrm{Ker}}
\newcommand{\Dim}{\mathrm{dim}}
\newcommand{\Codim}{\mathrm{codim}}
\newcommand{\ind}{\mathrm{ind}}

\newcommand{\Ddt}{\tfrac{\mathrm D}{\mathrm dt}}
\newcommand{\Ddtt}{\tfrac{\mathrm D^2}{\mathrm dt^2}}
\newcommand{\mul}{\mathrm{mul}}
\newcommand{\sgn}{\mathrm{sgn}}
\newcommand{\dd}{\mathrm d}
\newcommand{\Id}{\mathrm{Id}}

\newcommand{\Bils}{\mathrm B_{\mathrm{sym}}}
\newcommand{\gR}{g_{\scriptscriptstyle{\textrm R}}}

\documentclass[oneside,final,letterpaper,10pt]{amsart}

\usepackage{epsfig}
\usepackage{times}

\numberwithin{equation}{section}

\title[Bifurcation of semi-Riemannian Geodesics]{Spectral Flow, Maslov Index and\\ Bifurcation of
semi-Riemannian Geodesics}
\author[P.\ Piccione]{Paolo Piccione}
\address{Departamento de Matem\'atica,\hfill\break\indent  Universidade de S\~ao Paulo,
Brazil}
\email{piccione@ime.usp.br}
\urladdr{http://www.ime.usp.br/\~{}piccione}

\author[A. Portaluri]{Alessandro Portaluri}
\address{Dipartimento di Matematica,\hfill\break\indent Politecnico di Torino, Italy}
\email{portalur@calvino.polito.it}
\author[D. V. Tausk]{Daniel V. Tausk}
\address{Departamento de Matem\'atica,\hfill\break\indent  Universidade de S\~ao Paulo,
Brazil}
\email{tausk@ime.usp.br}
\urladdr{http://www.ime.usp.br/\~{}tausk}
\subjclass[2000]{58E10, 58J55, 53D12, 34K18, 47A53}

\thanks{The authors are grateful to Prof.\ J.\ Pejsachowicz for suggesting the problem
and for many stimulating discussions during this research project.}
\date{October 2002}

\begin{document}

% Theorems and such

\theoremstyle{plain}\newtheorem{teo}{Theorem}[section]
\theoremstyle{plain}\newtheorem{prop}[teo]{Proposition}
\theoremstyle{plain}\newtheorem{lem}[teo]{Lemma} 
\theoremstyle{plain}\newtheorem{cor}[teo]{Corollary} 
\theoremstyle{definition}\newtheorem{defin}[teo]{Definition} 
\theoremstyle{remark}\newtheorem{rem}[teo]{Remark} 
\theoremstyle{definition}\newtheorem{example}[teo]{Example}

\theoremstyle{plain}\newtheorem*{convention}{Convention}
\theoremstyle{definition}\newtheorem*{defin0}{Definition}

%%%%%

\begin{abstract}
We give a functional analytical proof of the equality
between the Maslov index of a semi-Riemannian geodesic
and the spectral flow of the path of self-adjoint
Fredholm operators obtained from the index form. This fact, together with 
recent results on the bifurcation for critical points of 
strongly indefinite functionals (see \cite{FitzPejsaRecht})
imply that each non degenerate and non null conjugate (or $P$-focal)
point along a semi-Riemannian geodesic is a bifurcation point.
%Other bifurcation problems in the context of semi-Riemannian
%geodesics are discussed.
\end{abstract}

\maketitle
\tableofcontents

%%%%%%%%%%%%%%%%
%%%%%%%%%%%%%%%%
\begin{section}{Introduction}\label{sec:intro}
 Let $(M,g)$ be a semi-Riemannian manifold and $p\in M$; a point $q\in M$ is conjugate
to $p$ if $q$ is a critical value of the exponential map $\exp_p$, i.e., if
the linearized geodesic map $\mathrm d\exp_p$ is not injective at $\exp_p^{-1}(q)$.
It is a natural question to ask whether the non injectivity at the linear level
implies non uniqueness of geodesics between two conjugate points. For instance,
two antipodal points on the Riemannian round sphere are joined by infinitely many geodesics; 
however, it is easy to produce examples of conjugate points in complete Riemannian
manifolds that are joined by a unique geodesic. 

In order to make a more precise sense of the above question,
first one has to observe that any
information obtained from the linearized geodesic equation can only be of {\em local\/} character,
which implies that one should not expect to detect the existence of a finite
number of geodesics between two  points along $\gamma$ 
by merely looking at the Jacobi equation. A similar situation occurs, for instance, 
when studying {\em cut points\/} along a Riemannian geodesic, that
are not necessarily related to conjugate points. On the other hand, in a number
of situations it is desirable to have a better picture of the
geodesic behavior near a conjugate point, and
in order to investigate this situation we introduce the notion of bifurcation point:
\begin{defin0}\label{thm:defbifpoint}
Let $(M,g)$ be a semi-Riemannian geodesic, $\gamma:[a,b]\to M$ be a geodesic
in $M$ and $t_0\in\left]a,b\right[$. The point $\gamma(t_0)$ is said to be a
{\em bifurcation point for $\gamma$\/} (see Figure~\ref{fig:bif1}) if there exists a sequence
$\gamma_n:[a,b]\to M$ of geodesics in $M$ and a sequence $(t_n)_{n\in\N}\subset\left]a,b\right[$
satisfying the following properties:
\begin{enumerate}
\item\label{itm:1defbifpt} $\gamma_n(a)=\gamma(a)$ for all $n$;
\item\label{itm:2defbifpt} $\gamma_n(t_n)=\gamma(t_n)$ for all $n$;
\item\label{itm:3defbifpt} $\gamma_n\to\gamma$ as $n\to\infty$;
\item\label{itm:4defbifpt} $t_n\to t_0$ (and thus $\gamma_n(t_n)\to\gamma(t_0)$) as $n\to\infty$.
\end{enumerate}
\end{defin0}
\begin{figure}
\begin{center}
\psfull
\epsfig{file=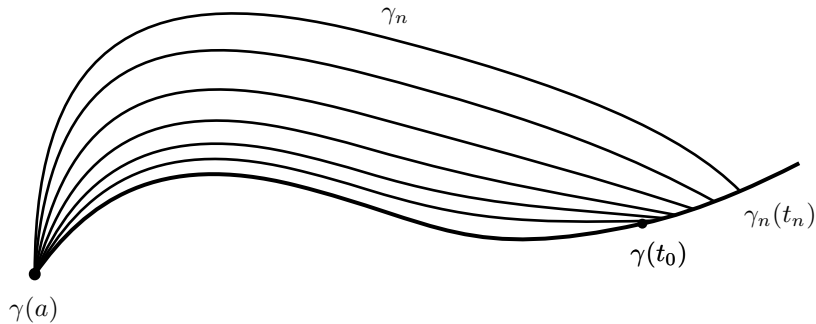} 
\caption{Bifurcation of geodesics.}\label{fig:bif1}
\end{center}
\end{figure}
The convergence of geodesics in condition \eqref{itm:3defbifpt} is meant in
any reasonable sense, for instance, it suffices to require that $\dot\gamma_n(a)\to\dot\gamma(a)$
as $n\to\infty$.

Using the Implicit Function Theorem, it follows immediately from the
above Definition that if $\gamma(t_0)$ is a bifurcation point for $\gamma$, then
necessarily $\gamma(t_0)$ must be conjugate to $\gamma(a)$ along $\gamma$. 
It is interesting to observe here that the above definition
of bifurcation point along a geodesic has strong analogies with 
Jacobi's original definition of conjugate point along an extremal of 
quadratic functionals (see for instance \cite[Definition~4, p.\ 114]{GelFom}).

The definition of
bifurcation point is well understood with the example of the 
paraboloid $z=x^2+y^2$, endowed with the
Euclidean metric of $\R^3$ (see Figure~\ref{fig:geoparaboloid}).
\begin{figure}
\begin{center}
\psfull
\epsfig{file=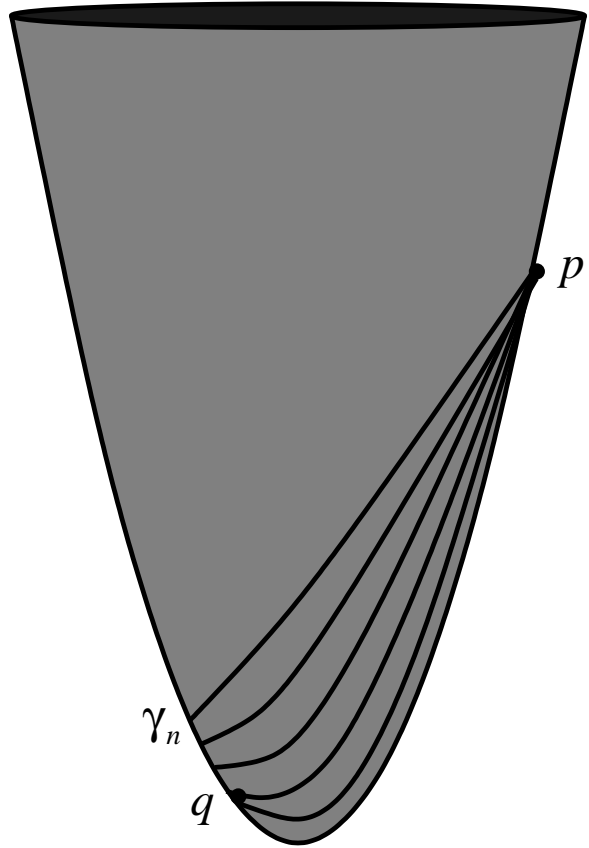} 
\caption{Geodesics issuing at a point $p$ of the paraboloid, tending to the
meridian through $p$.}\label{fig:geoparaboloid}
\end{center}
\end{figure}
Consider in this case the geodesic $\gamma$ given by the meridian issuing from a point 
$p$ distinct from the vertex of the paraboloid, with initial velocity pointing in the
negative $z$ direction. Such meridian goes downward towards the vertex, and then 
up again towards infinite on the opposite side of the paraboloid; this geodesic has a 
(unique) conjugate point $q$, and neighboring geodesics starting at $p$ intersect
the meridian at points $q_n\ne q$ that tend to $q$, and thus $q$ is a bifurcation 
point along $\gamma$.

Under the light of the above Definition, we reformulate the non uniqueness 
geodesic problem as follows: which conjugate points along a semi-Riemannian geodesic
are bifurcation points?  Several other bifurcation questions are naturally
associated to semi-Riemannian geometry. For instance, one could replace the
notion of conjugate point by that of {\em focal point\/} along a geodesic 
$\gamma$ relatively to an initial submanifold $P$ of $M$, and could ask which 
$P$-focal points are limits of endpoints  of geodesics starting orthogonally at $P$
and terminating on $\gamma$. 
%A somewhat different bifurcation problem consists in assuming that a given manifold
%$M$ is endowed with a one-parameter family $(g_\lambda)_{\lambda\in I}$ of semi-Riemannian
%metrics, and that it is given a curve $\gamma:[a,b]\to M$ which is a geodesic
%relatively to each $g_\lambda$. In this case, a natural question would be to
%determine whether for some $\lambda_*$ there exists a curve $\gamma_{\lambda_*}$
%joining $\gamma(a)$ and $\gamma(b)$ and distinct from $\gamma$ which is a geodesic
%relatively to $g_{\lambda_*}$. Again, if one wants to study this problem by
%means of the Jacobi equation, the correct question to ask
%is whether there exists a sequence $\lambda_n$ tending to $\lambda_*$ and a sequence
%$\gamma_n:[a,b]\to M$ of geodesics relatively to $g_{\lambda_n}$ joining
%$\gamma(a)$ and $\gamma(b)$, distinct
%from $\gamma$, and such that $\gamma_n\to\gamma$ as $n\to\infty$.

In this paper we use some recent results on  bifurcation theory
for strongly indefinite functionals (\cite{FitzPejsaRecht}) and on symplectic techniques 
for semi-Riemannian geodesics (\cite{CRASP, topology, fechado}) to give an answer to the
above questions. We outline briefly the ideas behind the theory of  Fitzpatrick, 
Pejsachowicz and Recht and how their result is employed in the present paper.
The most classical result on variational bifurcation (see \cite{Krasno}) states
that bifurcation for a smooth path of functionals having a trivial
branch of critical points with finite Morse index (assumed nondegenerate at the
endpoints) occurs at a given
singular critical point if such singular point determines a {\em jump\/}
of the Morse index. The variation  of the Morse index at the endpoints of a path
of essentially positive self-adjoint Fredholm operators is a homotopy invariant of the path;
recall to this aim that the space of essentially positive self-adjoint Fredholm 
operators form a contractible space, and that the invertible ones have an infinite
number of connected components, which are labelled by the Morse index. 
When dealing with strongly indefinite self-adjoint Fredholm operators,
then the topology of the space becomes richer (fundamental group isomorphic to $\Z$),
and no homotopy invariant for paths can be defined by simply looking at the
endpoints of the path. The {\em spectral flow\/} for a path, originally introduced
by  Atiyah, Patodi and Singer (see \cite{APS}), is an integer valued invariant
associated to paths of this type, and it is given, roughly speaking, 
by a signed count  of the eigenvalues that pass through zero at each
singular instants. The main result in \cite{FitzPejsaRecht} is that
bifurcation occurs at those singular instants whose contribution to the
spectral flow is non null (See Proposition~\ref{thm:FPR} below).

Consider now the geodesic bifurcation problem mentioned above.
By a suitable choice of coordinates in the space of
paths joining a fixed point $p$ in $M$ and a point variable along
a given geodesic $\gamma$ starting at $p$, the geodesic
bifurcation problem is reduced to a bifurcation problem
for a smooth family of strongly indefinite functionals defined 
in (an open neighborhood of $0$ of) a fixed Hilbert space.
The path of Fredholm operators corresponding to the index
form along the geodesic is studied, and the main result
of our computations is that its spectral flow coincides, up to a sign, with
another well known integer valued invariant of the geodesic,
called the Maslov index. Under a certain nondegeneracy assumption,
the Maslov index is computed as the sum of the signatures of all conjugate points
along the geodesic. Applying the theory of \cite{FitzPejsaRecht}, we get that nondegenerate conjugate
points with non vanishing signature are bifurcation points; more generally,
a bifurcation points is found in every segment of geodesic that contains a 
(possibly non discrete) set of conjugate points that give a non zero contribution 
to the Maslov index. In particular, Riemannian conjugate points are always
bifurcation points, as well as conjugate points along timelike or lightlike Lorentzian
geodesics. Similar results hold for focal points to an initial
nondegenerate submanifold. 
%As to the case of a one-parameter family of
%semi-Riemannian metrics $(g_\lambda)$, the condition of bifurcation
%at a value $\lambda_*$ for the geodesic $\gamma$ is that the
%$g_{\lambda_*}$-Maslov index of $\gamma$ be non zero.
\end{section}

\begin{section}{Fredholm bilinear forms on Hilbert spaces}
\label{sec:fredholmforms}
In this section we will discuss the notion of  index
of a Fredholm bilinear form on a Hilbert space relatively to
a closed subspace. The main goal (Proposition~\ref{thm:indrelindcoind}) is a result 
that gives the relative index of a form to the difference between
the index and the coindex
of suitable restrictions of the form.
\subsection{On the relative index of Fredholm forms}

Let $H$ be a Hilbert space with inner product $\langle\cdot,\cdot\rangle$,
and let $B$ a bounded symmetric bilinear form on $H$; there exists a unique
self-adjoint bounded operator $S:H\to H$ such that $B=\langle S\cdot,\cdot\rangle$,
that will be called the {\em realization of $B$\/} (with respect to $\langle\cdot,\cdot\rangle$).
$B$ is nondegenerate if its realization is injective, $B$ is strongly nondegenerate
if $S$ is an isomorphism. If $B$ is strongly nondegenerate, or if more generally 
$0$ is not an accumulation point of the spectrum of $S$, we will call the {\em negative
space\/} (resp., {\em the positive space}) of $B$  
the closed subspace $V^{\scriptscriptstyle-}(S)$ (resp., $V^{\scriptscriptstyle+}(S)$) of $H$ given by
$\chi_{\left]-\infty,0\right[}(S)$ (resp., $\chi_{\left]0,+\infty\right[}(S)$), where $\chi_I$ denotes
the characteristic function of the interval $I$.
We will say that $B$ is {\em Fredholm\/} if $S$ is Fredholm, or that $B$ is {\em RCPPI},
{\em realized by a compact perturbation of a positive isomorphism}, (resp., {\em RCPNI})
if $S$ is of the form $S=P+K$ (resp., $S=N+K$) where $P$ is a positive isomorphism
of $H$ ($N$ is a negative isomorphism of $H$) and $K$ is compact. 
Observe that
the properties of being Fredholm, RCPPI or RCPNI do not depend on the inner product, although
the realization $S$ and the spaces $V^\pm(S)$ do.

The index (resp., the coindex) of $B$, denoted by $n_-(B)$ (resp., $n_+(B)$) is the dimension of
$V^{\scriptscriptstyle-}(S)$ (resp., of
$V^{\scriptscriptstyle+}(S)$); the {\em nullity\/} of $B$, denoted by $n_0(B)$ is the dimension of the
kernel of $S$.

If $B$ is RCPPI (resp., RCPNI), then both its nullity $n_0(B)$ and its index $n_-(B)$ (resp., and its
coindex $n_+(B)$) are finite numbers. 

Given a closed subspace $W\subset H$, the {\em $B$-orthogonal complement of $W$}, denoted by
$W^{\perp_B}$, is the closed subspace of $H$:
\[W^{\perp_B}=\big\{x\in H:B(x,y)=0\ \text{for all}\ y\in W\big\};\]
clearly,
\[W^{\perp_B}=S^{-1}(W^\perp).\]
If $B$ is  Fredholm, $S$ is its realization and $W\subset H$ is any subspace, then the following
properties hold:
\begin{enumerate}
\item $B$ is nondegenerate iff it is strongly nondegenerate;
\item $n_0(B)<+\infty$;
\item $(W^{\perp_B})^{\perp_B}=\overline W+\mathrm{Ker}(S)$;
\item if $W$ is closed, then $W+W^{\perp_B}$ is closed;
\item if $W$ is closed and $B\vert_W$ (i.e., the restriction of $B$ to $W\times W$)
in nondegenerate, then also $B\vert_{W^{\perp_B}}$ is nondegenerate and $H=W\oplus W^{\perp_B}$.
\end{enumerate}

Let us now recall a few basic things on the notion of commensurability of
closed subspaces (see reference \cite{Abbo1} for more details). 
Let $V,W\subset H$ be closed subspaces and let $P_V$ and
$P_W$ denote the orthogonal projections respectively onto $V$ and $W$. We say that
$V$ and $W$ are {\em commensurable\/} if $P_V-P_W$ is a compact operator.
Equivalently, $V$ and $W$ are commensurable if both $P_{W^\perp}P_V$ and $P_{V^\perp}P_W$
are compact; if $V$ and $W$ are commensurable the {\em relative dimension\/}
$\mathrm {dim}_V(W)$ of $W$ with respect to $V$ is defined as:
\[\mathrm{dim}_V(W)=\mathrm{dim}(W\cap V^\perp)-\mathrm{dim}(W^\perp\cap V).\]
Clearly, if $V$ and $W$ are commensurable, then $V^\perp$ and $W^\perp$ are commensurable, and:
\[\mathrm{dim}_{V^\perp}(W^\perp)=-\mathrm{dim}_V(W).\]
The notion of commensurability of subspaces does not depend on the Hilbert space inner product
of $H$.

\begin{prop}
\label{thm:era1.6}
Let $S,T$ be linear bounded self-adjoint operators on $H$ whose difference
$K=S-T$ is compact. Then $V^{\scriptscriptstyle-}(S)$ (resp., $V^{\scriptscriptstyle+}(S)$)
is commensurable with $V^{\scriptscriptstyle-}(T)$ (resp., with $V^{\scriptscriptstyle+}(T)$).

Conversely, assume that $S$ is a bounded self-adjoint Fredholm operator on $H$, and let
$H=W^{\scriptscriptstyle-}\oplus W^{\scriptscriptstyle+}$ be an orthogonal decomposition 
of $H$ such that $W^{\scriptscriptstyle-}$ is commensurable with $V^{\scriptscriptstyle-}(S)$
and $W^{\scriptscriptstyle+}$ is commensurable with $V^{\scriptscriptstyle+}(S)$.
Then there exists an invertible self-adjoint operator $T$ on $H$ such that
$V^{\scriptscriptstyle-}(T)=W^{\scriptscriptstyle-}$,
$V^{\scriptscriptstyle+}(T)=W^{\scriptscriptstyle+}$ and such that $S-T$ is compact. 
\end{prop}
\begin{proof}
See \cite[Proposition~2.3.2 and Proposition~2.3.5]{Abbo1}.
\end{proof}

\begin{lem}\label{thm:RCP-comm}
Let $B$ be a Fredholm symmetric bilinear form on the Hilbert space $H$
and let  $W\subset H$ be a closed
subspace. 
Then, the following are equivalent:
\begin{itemize}
\item[(a)]  $B\vert_W$ is RCPNI and $B\vert_{W^{\perp_B}}$ is RCPPI;
\item[(b)] there exists a Hilbert space inner product $\langle\cdot,\cdot\rangle$ on $H$ 
such that $W$ is commensurable with $V^{\scriptscriptstyle-}(S)$, where $S$ is the realization of $B$
with  respect to
$\langle\cdot,\cdot\rangle$.
\end{itemize}
\end{lem}
\begin{proof}
Assume that (b) holds; fix a Hilbert space inner product $\langle\cdot,\cdot\rangle$
in $H$ and let $S$ be the realization of $B$ with respect to $\langle\cdot,\cdot\rangle$
so that $W$ is commensurable with $V^{\scriptscriptstyle-}(S)$. Then $W^\perp$  is commensurable with
$V^{\scriptscriptstyle-}(S)^\perp= V^{\scriptscriptstyle+}(S)\oplus\Ker(B)$. Moreover, since
$\Ker(B)$ is finite dimensional, then $W^\perp$ is also commensurable with
$V^{\scriptscriptstyle+}(S)$. By Proposition~\ref{thm:era1.6}, there exists an invertible self-adjoint
operator $T:H\to H$ such that $V^{\scriptscriptstyle-}(T)=W$, $V^{\scriptscriptstyle+}(T)=W^\perp$,
and with $S=T+K$, with $K$ compact. It follows easily that $B\vert_W$ is RCPNI
(namely, if $P$ denotes the orthogonal projection onto $W$, the realization
of $B\vert_W$ is $PS\vert_W=(PT+PK)\vert_W=(T+PK)\vert_W$), and $B\vert_{W^\perp}$ is RCPPI.
Observe in particular that $W\cap W^{\perp_B}=\Ker(B\vert_W)$ is finite dimensional.
To prove that $B\vert_{W^{\perp_B}}$ is RCPPI we argue as follows; denote by $P$
the orthogonal projection onto $W$ and by $P^{\scriptscriptstyle\perp}=1-P$ the orthogonal
projection onto $W^\perp$. As we
have observed, $W^{\perp_B}=S^{-1}(W^\perp)$; hence,
for all $x,y\in W^{\perp_B}$ we have:
\begin{equation}\label{eq:era6}\begin{split}B(x,y)&\,=\langle Sx,y\rangle=\langle
Sx,P^{\scriptscriptstyle\perp} y\rangle=\langle SPx,P^{\scriptscriptstyle\perp} y\rangle+\langle 
SP^{\scriptscriptstyle\perp} x,P^{\scriptscriptstyle\perp} y\rangle=\\ &\,=\langle
P^{\scriptscriptstyle\perp}KPx,y\rangle +\langle
P^{\scriptscriptstyle\perp}TP^{\scriptscriptstyle\perp}x,y\rangle+
\langle P^{\scriptscriptstyle\perp}KP^{\scriptscriptstyle\perp}x,y\rangle.\end{split}
\end{equation}
In the above equality we have used the fact that $W$ and $W^\perp$ are $T$-invariant.
From \eqref{eq:era6} we deduce that $B\vert_{W^{\perp_B}}$ is represented by
a compact perturbation of the operator
$\widetilde T:W^{\perp_B}\to W^{\perp_B}$ given by $\widetilde
T=P^{\scriptscriptstyle\perp_B}P^{\scriptscriptstyle\perp}TP^{\scriptscriptstyle\perp}
\vert_{W^{\perp_B}}$ (where $P^{\scriptscriptstyle\perp_B}$ is the orthogonal
projection onto $W^{\perp_B} $) which is positive
semi-definite. The kernel of $\widetilde T$ is easily computed as the finite
dimensional space $W^{\perp_B}\cap T^{-1}\big(W\cap W^{\perp_B}\big)$; it follows that $\widetilde
T$ is a compact perturbation of a positive isomorphism of $W^{\perp_B}$,
which proves that (b) implies (a).

Conversely, if $B\vert_W$ is RCPNI and $B\vert_{W^{\perp_B}}$ is RCPPI, then
clearly $W_1=W\cap W^{\perp_B}$ is finite dimensional; let $\widetilde W$ be
any closed complement of $W_1$ in $W$. It follows that $B\vert_{\widetilde W}$ is
nondegenerate, which implies that we have a direct sum decomposition
$H=\widetilde W\oplus\widetilde W^{\perp_B}$. If $\langle\cdot,\cdot\rangle$
is any Hilbert space inner product for which $\widetilde W$ and $\widetilde W^{\perp_B}$
are orthogonal, then it is easily checked that the corresponding realization
$S$ of $B$ is such that $V^{\scriptscriptstyle-}(S)$ is commensurable with $W$.

This concludes the proof.
\end{proof}

Assume now that $B$ is a symmetric bilinear form, $S$ is its realization; if $W$ is closed
subspace of $H$ which is commensurable with $V^{\scriptscriptstyle-}(S)$, the one defines the {\em
relative index\/} of $B$ with respect to $W$, denoted by $\ind_W(B)$, the integer number:
\[\ind_W(B)=\mathrm{dim}_W\big(V^{\scriptscriptstyle-}(S)\big).\]
Again, the relative index is independent of the inner product, and the
following equality holds:
\[\ind_W(B)=\sup\big\{\mathrm{dim}_W(V): V\ \text{is commensurable with}\ 
V^{\scriptscriptstyle-}(S)\big\}.\]

\subsection{Computation of the relative index}
A subspace $Z$ of $H$ is said to be {\em isotropic\/} for the symmetric bilinear form $B$ of
$B\vert_Z\equiv0$.
\begin{lem}\label{thm:dimisotropic}
Let $B$ be a RCPPI symmetric bilinear form on $H$, and let $Z\subset H$ be an isotropic subspace of
$B$. Then:
\[n_-(B)=n_-\big(B\vert_{Z^{\perp_B}}\big)+\mathrm{dim}(Z).\]
\end{lem}
\begin{proof}
Since $B$ is RCPPI, then the index $n_-(B)$ is finite, and so 
$n_-\big(B\vert_{Z^{\perp_B}}\big)$ and $\Dim(Z)$ are finite.
Clearly, $Z\subset Z^{\perp_B}$; let $U\subset Z^{\perp_B}$ be a closed subspace 
such that $Z^{\perp_B}=Z\oplus U$, so that $B\vert_U$ is nondegenerate and
$H=U\oplus U^{\perp_B}$. Moreover:
\[n_-(B)=n_-\big(B\vert_U\big)+n_-\big(B\vert_{U^{\perp_B}}\big).\]
Since $Z$ is isotropic, then $n_-\big(B\vert_U\big)=n_-\big(B\vert_{Z^{\perp_B}}\big)$;
to conclude the proof we need to show that $n_-\big(B\vert_{U^{\perp_B}}\big)=
\Dim(Z)$. To this aim, observe first that $\Dim(U^{\perp_B})=2\Dim(Z)$.
Namely, $\Dim(U^{\perp_B})=\Codim(U)$; moreover, $\Codim_{Z^{\perp_B}}(U)=
\Dim(Z)$, and $\Codim(Z^{\perp_B})=\Dim(Z)$.
Thus, keeping in mind that the dimension of an isotropic subspace is
less than or equal to the index and the coindex, we have:
\[n_-\big(B\vert_{U^{\perp_B}}\big)+n_+\big(B\vert_{U^{\perp_B}}\big)=\Dim(
U^{\perp_B})=2\,\Dim(Z)\le n_-\big(B\vert_{U^{\perp_B}}\big)+n_+\big(B\vert_{U^{\perp_B}}\big),
\]
which proves that $n_-\big(B\vert_{U^{\perp_B}}\big)=n_+\big(B\vert_{U^{\perp_B}}\big)=
\Dim(Z)$ and concludes the proof.
\end{proof}

\begin{lem}\label{thm:indtilde}
Let $B$ be a nondegenerate Fredholm symmetric bilinear form on $H$ and $W\subset H$ be a closed
subspace
such that $B\vert_{W^{\perp_B}}$ is RCPPI. Let $\widetilde W$ be any closed complement\footnote{%
for instance, $\widetilde W$ is the orthogonal complement of $W\cap W^{\perp_B}$ in $W$
with respect to any inner product.}
of $W\cap W^{\perp_B}$ in $W$. Then the following identity holds:
\[n_-\big(B\vert_{\widetilde W^{\perp_B}}\big)=n_-\big(B\vert_{ 
W^{\perp_B}}\big)+\mathrm{dim}\big(W\cap W^{\perp_B}\big).\]
\end{lem}
\begin{proof}
We start with the observation that $\Ker\big(B\vert_{W}\big)=\Ker\big(B\vert_{W^{\perp_B}}\big)=W\cap
W^{\perp_B}$; this implies in particular that $B\vert_{\widetilde W}$  and $B\vert_{\widetilde
W^{\perp_B}}$ are nondegenerate. 
Since $B\vert_{W^{\perp_B}}$ is RCPPI, then $n_-\big(B\vert_{W^{\perp_B}}\big)
$ and $\Dim(W\cap W^{\perp_B})=n$ are finite numbers.

Since $\Codim_{\widetilde W^{\perp_B}}\big(W^{\perp_B}\big)=n$, then:
\[n_-\big(B\vert_{\widetilde W^{\perp_B}}\big)\le n_-\big(B\vert_{W^{\perp_B}}\big)+n,\]
from which it follows that $ n_-\big(B\vert_{\widetilde W^{\perp_B}}\big)$ is finite;
moreover,
$B\vert_{\widetilde W^{\perp_B}}$ is RCPPI. The conclusion now follows
easily from
Lemma~\ref{thm:dimisotropic}, applied to the bilinear form $B\vert_{\widetilde W^{\perp_B}}$ and the
isotropic space $Z=W\cap W^{\perp_B}$.
\end{proof}

We are finally ready to give our central result concerning
the computation of the relative index of a Fredholm bilinear
form $B$  in terms of index and coindex of suitable restrictions
of $B$:

\begin{prop}\label{thm:indrelindcoind}
Let $B$ be a Fredholm symmetric bilinear form on $H$, $S$ its realization and let $W\subset H$ be
a closed subspace which is commensurable with $V^{\scriptscriptstyle-}(S)$. Then the relative index
$\ind_W(B)$ is given by:
\begin{equation}\label{eq:era7}
\ind_W(B)=n_-\big(B\vert_{W^{\perp_B}}\big)-n_+\big(B\vert_W\big).
\end{equation}
\end{prop}
\begin{proof}
Assume first that $B$ is nondegenerate on $W$; then
have a direct sum decomposition $H=W\oplus W^{\perp_B}$. The relative
$\ind_{W}(B)$ does not change if we change the inner product of
$H$; we can therefore assume that $W$ and $W^{\perp_B}$ are orthogonal
subspaces of $H$. Then, $S=S^{\scriptscriptstyle-}\oplus S^{\scriptscriptstyle+}$,
where $S^{\scriptscriptstyle-}:W\to W$ is
the realization of $B\vert_{W}$
and $S^{\scriptscriptstyle+}:W^{\perp_B}\to W^{\perp_B}$ is the realization of
$B\vert_{W^{\perp_B}}$. Moreover,
$V^{\scriptscriptstyle-}(S)=V^{\scriptscriptstyle-}(S^{\scriptscriptstyle-})\oplus
V^{\scriptscriptstyle-}(S^{\scriptscriptstyle+})$. An immediate calculation  yields:
\[\begin{split} \ind_{W}(B)\,&=
\Dim\big(V^{\scriptscriptstyle-}(S)\cap W^{\perp_B}\big)-\Dim\big(
V^{\scriptscriptstyle-}(S)^\perp\cap W\big)
\\ &=\Dim\big(V^{\scriptscriptstyle-}(S)\cap
W^{\perp_B}\big)-\Codim_{W}\big(V^{\scriptscriptstyle-}(S^{\scriptscriptstyle-})\big) \\
&=\Dim\big(V^{\scriptscriptstyle-}(S^{\scriptscriptstyle+})\big)
-\Codim_{W}\big(V^{\scriptscriptstyle-}(S^{\scriptscriptstyle-})\big) \\
&=n_-\big(B\vert_{W^{\perp_B}}\big)-n_+\big(B\vert_{W}\big).\end{split}\]

Let us consider now the case that $B\vert_W$ is degenerate; by Lemma~\ref{thm:RCP-comm},
$B\vert_W$ is RCPNI, and so  $\Dim\big(W\cap W^{\perp_B}\big)=n<+\infty$.
Set $\widetilde W=\big(W\cap W^{\perp_B}\big)^\perp\cap W$, so that $B\vert_{\widetilde W}$
is nondegenerate; moreover, $V^{\scriptscriptstyle-}(S)$ is commensurable
with $\widetilde W$, because it has finite codimension in $W$. 
We can then apply the first part of the proof, and we obtain:
\begin{equation}\label{eq:era8}\ind_{\widetilde W}(B)=n_-\big(B\vert_{\widetilde
W^{\perp_B}}\big)-n_+\big(B\vert_{
\widetilde W}\big).\end{equation}
Clearly, 
\begin{equation}\label{eq:era9}
n_+\big(B\vert_{
\widetilde W}\big)=n_+\big(B\vert_{W}\big);\end{equation}
moreover, by definition of relative index:
\begin{equation}\label{eq:era10}
\ind_{\widetilde W}(B)=\ind_{W}(B)+n.\end{equation}
Finally, by Lemma~\ref{thm:RCP-comm}, $B\vert_{W^{\perp_B}}$ is RCPPI, and by Lemma~\ref{thm:indtilde}:
\begin{equation}\label{eq:era11}n_-\big(B\vert_{\widetilde
W^{\perp_B}}\big)=n_-\big(B\vert_{W^{\perp_B}}\big)+n.\end{equation} 
Formulas \eqref{eq:era8}, \eqref{eq:era9}, \eqref{eq:era10} and \eqref{eq:era11}
yield \eqref{eq:era7} and conclude the proof.
\end{proof}

\end{section}

%%%%%%%%%%%%%%%%
%%%%%%%%%%%%%%%%
\begin{section}{On the spectral flow of a path of self-adjoint Fredholm operators}
\label{sec:spectralflow}

In this section we will recall some facts from the theory
of variational bifurcation for strongly indefinite 
functionals. The basic reference for the material
presented is \cite{FitzPejsaRecht}; as to the definition
and the basic properties of the spectral flow we refer
to the nice article by Phillips \cite{Phillips}, from which we will
borrow some of the notations.
\subsection{Spectral flow}
Let us consider an infinite dimensional
 separable real Hilbert space $H$. We will denote by $\Bcal(H)$ and
$\Kcal(H)$ respectively the algebra of all bounded linear operators
on $H$ and the closed two-sided ideal of $\Bcal(H)$ consisting 
of all compact operators on $H$;
the Calkin algebra $\Bcal(H)/\Kcal(H)$ will be denoted by $\Qcal(H)$, and
$\pi:\Bcal(H)\to\Qcal(H)$ will denote the quotient map.
The {\em essential spectrum\/} $\sigma_\ess(T)$ of a bounded linear operator $T\in\Bcal(H)$
is the spectrum of $\pi(T)$ in the Calkin algebra $\Qcal(H)$.
Let $\Fcal(H)$ and $\Fcalsa(H)$ denote respectively the space
of all Fredholm (bounded) linear operators on $H$ and the space
of all self-adjoint ones. An element $T\in\Fcalsa(H)$ is
said to be {\em essentially positive} (resp., {\em essentially negative})
if $\sigma_\ess(T)\subset\R^+$ (resp., if $\sigma_\ess(T)\subset\R^-$),
and {\em strongly indefinite\/} if it is neither essentially
positive nor essentially negative.

The symbols $\Fcalsa_+(H)$, $\Fcalsa_-(H)$ and $\Fcalsa_*(H)$ will denote  
the subsets  of $\Fcalsa(H)$ consisting respectively of all essentially
positive, essentially negative and strongly indefinite self-adjoint
Fredholm operators on $H$.
These sets are precisely the three connected components
of $\Fcalsa(H)$; $\Fcalsa_+(H)$ and $\Fcalsa_-(H)$ are contractible,
while $\Fcalsa_*(H)$ is homotopically equivalent to $U(\infty)=\lim_nU(n)$,
and it has infinite cyclic fundamental group.

Given a continuous path $S:[0,1]\to\Fcalsa_*(H)$ with $S(0)$ and $S(1)$ 
invertible, the {\em spectral flow\/} of $S$, denoted by $\spfl(S)$, is an integer number
which is given, roughly speaking, by the net number of eigenvalues that pass through zero in the
positive direction from the start of the path to its end. There exist  several
equivalent definitions of the spectral flow in the literature; we like to mention here
the definition given in \cite{Phillips} using functional calculus, and that reduces
the problem to a simple dimension counting of finite
rank projections.

More precisely, let $\chi_I$ denote the characteristic function
of the interval $I$; for all $S\in\Fcalsa_*(H)$ there exists $a>0$ and
a neighborhood $U$ of $S$ in $\Fcalsa_*(H)$ such that the map
$T\mapsto\chi_{[-a,a]}(T)$ is norm continuous in $U$, and it takes values
in the set of projections of {\em finite\/} rank.
Denote by  $C^0_\#\big([0,1],\Fcalsa_*(H)\big)$ the set of all continuous paths
$S:[0,1]\to \Fcalsa_*(H)$ such that $S(0)$ and $S(1)$ are invertible. Given
$S\in C^0_\#\big([0,1],\Fcalsa_*(H)\big)$, then by the above property
one can choose a partition $0=t_0<t_1<\ldots<t_N=1$ of $[0,1]$
and positive numbers $a_1,\ldots,a_N$
such that the maps $t\mapsto\chi_{[-a_i,a_i]}\big(S(t)\big)$ are
continuous and of finite rank on $[t_{i-1},t_i]$ for all $i$.
The spectral flow of the path $S$ is defined to be the sum:
\[\sum_{i=1}^n\Big[\mathrm{rk}\big(\chi_{[0,a_i]}(S(t_i))\big)-\mathrm{rk}
\big(\chi_{[0,a_i]}(S(t_{i-1}))\big)\Big],\]
where $\mathrm{rk}$ is the rank of a projection.
With the above formula, the spectral flow is well defined, i.e., it does
not depend on the choice of the partition $(t_i)$ and
of the positive numbers $(a_i)$, and the map $\spfl:C^0_\#\big([0,1],\Fcalsa_*(H)\big)\to
\Z$ has the following properties:
\begin{itemize}
\item it is additive by concatenation;
\item if $S\in C^0_\#\big([0,1],\Fcalsa_*(H)\big)$ is such that $S(t)$ is
invertible for all $t$, then $\spfl(S)=0$;
\item it is invariant by homotopies with fixed endpoints;
\item the induced map $\spfl:\pi_1\big(\Fcalsa_*(H)\big)\to\Z$ is an isomorphism.
\end{itemize}

For the purposes of the present paper, it will be useful
to give a different description of the spectral flow,
which follows the approach in \cite{FitzPejsaRecht}.
As we have observed, $\Fcalsa_*(H)$ is not simply connected,
and therefore no non trivial homotopic invariant for
curves in $\Fcalsa_*(H)$ can be defined only
in terms of the value at the endpoints. However,
in \cite{FitzPejsaRecht} it is shown that the spectral
flow can be defined in terms of the endpoints, provided that 
the path $S$ has the special form $S(t)=\mathfrak J+K(t)$, where
$\mathfrak J$ is a fixed symmetry of $H$ and $t\mapsto K(t)$ is a path of
compact operators.
By a {\em symmetry\/} of the Hilbert space $H$ it is meant 
an operator $\mathfrak J$ of the form \[\mathfrak J=P_+-P_-,\]
where $P_+$ and $P_-$ are the orthogonal projections
onto infinite dimensional closed subspaces $H_+$ and $H_-$ of $H$
such that $H=H_+\oplus H_-$; assume that such a symmetry $\mathfrak J$ has
been fixed.

Denote by $\Bcal_o(H)$ the group of all invertible elements
of $\Bcal(H)$. There is an action of $\Bcal_o(H)$ on $\Fcalsa(H)$
given by:
\[\Bcal_o(H)\times \Fcalsa(H)\ni (M,S)\longmapsto M^*SM\in \Fcalsa(H); \]
this action preserves the three connected components of $\Fcalsa(H)$.
Two elements in the same orbit are said to be {\em cogredient};
the orbit of each element in
$\Fcalsa_*(H)$ meets the affine space
$\mathfrak J+\Kcal(H)$, i.e., given any $S\in\Fcalsa_*(H)$ there exists
$M\in\Bcal_o(H)$ such that $M^*SM=\mathfrak J+K$, where $K$ is compact.
Moreover, using a suitable fiber bundle structure and standard lifting arguments,
it is shown in \cite{FitzPejsaRecht} that if $t\mapsto S(t)\in\Fcalsa_*(H)$ is a
path of class $C^k$, $k=0,\ldots,+\infty$, then one can find a $C^k$ curve
$t\mapsto M(t)\in\Bcal_o(H)$ such that $M(t)^*S(t)M(t)=\mathfrak J+K(t)$,
where $t\mapsto K(t)$ is a $C^k$ curve  of compact operators.
Among the central results of \cite{FitzPejsaRecht} the authors prove
that the spectral flow of a path of strongly indefinite self-adjoint
Fredholm operators is invariant by cogredience, and that for
paths that are compact perturbation of a fixed symmetry the
spectral flow is given as the relative dimension of the
negative eigenspaces at the endpoints:
\begin{prop}\label{thm:cogrinvformsf}
Let $S:[0,1]\to\Fcalsa_*(H)$  be a continuous path such that 
$S(0)$ and $S(1)$ are invertible, denote by $B(t)=\langle S(t)\cdot,\cdot\rangle$
the corresponding bilinear form on $H$, and let
$M:[0,1]\to\Bcal_o(H)$ be a continuous curve with $L(t):=M(t)^*S(t)M(t)$ of the form
$\mathfrak J+K(t)$, with $K(t)$ compact for all $t$. Then:
\begin{enumerate}
\item $\spfl(S)=\spfl(L)$;
\item\label{itm:2calcspfl} $\spfl(L)= \ind_{V^{\scriptscriptstyle-}
\big(L(1)\big)}\big(B(0)\big)\\ \phantom{\spfl(L)}=\Dim\Big(V^{\scriptscriptstyle-}\big(L(0)\big)\cap
V^{\scriptscriptstyle+}
\big(L(1)\big)\Big)-\Dim\Big(V^{\scriptscriptstyle+}\big(L(0)\big)\cap V^{\scriptscriptstyle-}
\big(L(1)\big)\Big)$.
\end{enumerate}
\end{prop}
\begin{proof}
See \cite[Proposition~3.2, Proposition~3.3]{FitzPejsaRecht}.
\end{proof}
Observe that, since $\Dim_W(V)=-\Dim_V(W)$, the equality in part~\eqref{itm:2calcspfl}
of Proposition~\ref{thm:cogrinvformsf} can be rewritten as:
\begin{equation}\label{thm:cogrinvformsf2}
\spfl(L)=-\ind_{V^{\scriptscriptstyle-}
\big(L(0)\big)}\big(B(1)\big)
\end{equation}

\subsection{Bifurcation for a path of strongly indefinite functionals}
Let $H$ be a real separable Hilbert space, $U\subset H$ a neighborhood of
$0$ and $f_\lambda:U\to\R$ a family of smooth (i.e., of class $C^2$) functionals 
depending smoothly on $\lambda\in[0,1]$.
Assume that $0$ is a critical point of $f_\lambda$ for all $\lambda\in[0,1]$.
An element $\lambda_*\in[0,1]$ is said to be a {\em bifurcation value\/}
if there exists a sequence $(\lambda_n)_n$ in $[0,1]$ and a sequence
$(x_n)_n\in U$ such that:
\smallskip

\begin{enumerate}
\item $x_n$ is a critical point of $f_{\lambda_n}$ for all $n$;
\smallskip

\item $x_n\ne0$ for all $n$ and $\lim\limits_{n\to\infty}x_n=0$;
\smallskip

\item $\lim\limits_{n\to\infty}\lambda_n=\lambda_*$.
\end{enumerate}
The main result concerning the existence of a bifurcation
value for a path of strongly indefinite functionals is
the following:
\begin{prop}\label{thm:FPR}
Let $S(\lambda)=\mathrm d^2f_\lambda(0)$ be the continuous path of
self-adjoint Fredholm operators on $H$ given by the second variation
of $f_\lambda$ at $0$. Assume that $S$ takes values in $\Fcalsa_*(H)$
for all $\lambda\in[0,1]$, and that $S(0)$ and $S(1)$ are invertible.
If $\spfl(S)\ne0$, then there exists a bifurcation value $\lambda_*\in\left]0,1\right[$.
\end{prop}
\begin{proof}
See \cite[Theorem~1]{FitzPejsaRecht}.
\end{proof}
It is  obvious that, being a local notion, bifurcation
can be defined also in the case of a smooth family of $C^2$-functionals
$f_\lambda$, $\lambda\in[a,b]$, defined on (an open subset of) a Hilbert
manifold $\Omega$, in the case that there exists a common critical
point $\mathfrak z\in\Omega$ for all the $f_\lambda$'s.
Using local charts around $\mathfrak z$ (and thus identifying 
the tangent spaces at each point near $\mathfrak z$ with a fixed Hilbert
space) one sees immediately that the result of Proposition~\ref{thm:FPR}
holds also in this setting.  On the other hand, global existence
results for nontrivial branches of critical points in the linear case
cannot be extended directly to the case of manifolds.
\end{section}

%%%%%%%%%%%%%%%%
%%%%%%%%%%%%%%%%
\begin{section}{On the Maslov index}
\label{sec:maslov}

We will henceforth consider a smooth manifold $M$ endowed with a semi-Riemannian
metric tensor $g$; by the symbol $\Ddt$ we will denote the covariant differentiation
of vector fields along a curve in the Levi--Civita connection of $g$, while
$R$ will denote the curvature tensor of this connection chosen with the sign convention: $R(X,Y)=
[\nabla_X,\nabla_Y]-\nabla_{[X,Y]}$. Set $n=\Dim(M)$.
\subsection{Semi-Riemannian conjugate points}
\label{sub:conjugatepts}
Let $\gamma:[0,1]\to M$ be a geodesic in $(M,g)$; consider the Jacobi equation 
for vector fields  along $\gamma$:
\begin{equation}\label{eq:Jacobieq}
\Ddtt J-R(\dot\gamma,J)\,\dot\gamma=0.
\end{equation}
Let $\mathbb J$ denote the $n$-dimensional space:
\begin{equation}\label{eq:defJbb}
\mathbb J=\big\{J\ \text{solution of \eqref{eq:Jacobieq} such that}\ 
J(0)=0\big\}.\end{equation}
 A point $\gamma(t_0)$, $t_0\in\left]0,1\right]$ is said to be {\em conjugate\/}
to
$\gamma(0)$ if there exists a non zero   $J\in\mathbb J$ such that $J(t_0)=0$.

Set $\mathbb J[t_0]=\big\{J(t_0):J\in\mathbb J\big\}$; the codimension of
$\mathbb J[t_0]$ in $T_{\gamma(t_0)}M$ is called the {\em multiplicity\/}
of the conjugate point $\gamma(t_0)$, denoted by $\mul(t_0)$. The signature of
the restriction of $g$ to the $g$-orthogonal complement $\mathbb J[t_0]^\perp$
is called the {\em signature\/} of $\gamma(t_0)$, and will be denoted by $\sgn(t_0)$.
The conjugate point $\gamma(t_0)$ is said to be {\em nondegenerate\/} if such restriction
is nondegenerate; clearly, if $g$ is Riemannian (i.e., positive definite) then
every conjugate point is nondegenerate and its signature coincides with its multiplicity 
(the same is true for conjugate points along timelike or lightlike Lorentzian geodesics,
see the proof of Corollary~\ref{thm:corRiemcausLor}).

It is well known that nondegenerate conjugate points are isolated, while the
distribution of degenerate conjugate points can be quite arbitrary (see \cite{fechado}).
\subsection{The Maslov index: geometrical definition.}
\label{sub:maslovgeom}
Let $v_1,\ldots,v_n$ be a $g$-orthonormal basis of $T_{\gamma(0)}M$
and consider the parallel frame $V_1,\ldots,V_n$  obtained by parallel
transport of the $v_i$'s along $\gamma$. This frame gives us
isomorphisms $T_{\gamma(t)}M\to\R^n$ that carry the metric tensor $g$
to a {\em fixed\/} symmetric bilinear form on $\R^n$, still denoted by $g$.
Observe that, by the choice of a parallel trivialization of the tangent bundle
$TM$ along $\gamma$, covariant differentiation for vector fields along
$\gamma$ corresponds to standard differentiation of $\R^n$-valued maps,
and the Jacobi equation  \eqref{eq:Jacobieq} becomes the Morse--Sturm system:
\begin{equation}\label{eq:MS}
J''=RJ,
\end{equation}
where $R$ is a smooth curve of $g$-linear endomorphisms of $\R^n$.

%%%%%%
%%%%%%
%%%%%%
Consider the space $\R^n\oplus{\R^n}^*$ endowed with the {\em canonical symplectic form\/}
\[\omega\big((v_1,\alpha_1),(v_2,\alpha_2)\big)=\alpha_2(v_1)-\alpha_1(v_2),\quad v_1,v_2\in\R^n,\
\alpha_1,\alpha_2\in{\R^n}^*.\]
We denote by $\Spl(2n,\R)$ the {\em symplectic group\/} of $\R^n\oplus{\R^n}^*$, i.e., the Lie group of
all symplectomorphisms of $\R^n\oplus{\R^n}^*$; by $\spl(2n,\R)$ we denote the {\em Lie algebra\/} of
$\Spl(2n,\R)$. Recall that a {\em Lagrangian\/} subspace
$L$ of
$\R^n\oplus{\R^n}^*$ is an
$n$-dimensional subspace on which $\omega$ vanishes. We denote by $\Lambda$ the 
{\em Lagrangian Grassmannian\/}  of $\R^n\oplus{\R^n}^*$ which is the set of all
Lagrangian subspaces of $\R^n\oplus{\R^n}^*$. The Lagrangian Grassmannian is a
$\frac12n(n+1)$-dimensional compact and connected real-analytic embedded submanifold of the
Grassmannian of all $n$-dimensional subspaces of $\R^n\oplus{\R^n}^*$. Given a Morse--Sturm
system
\eqref{eq:MS} we set:
\begin{equation}\label{eq:defellt}
\ell(t)=\big\{\big(J(t),gJ'(t)\big):J\in\mathbb J\big\}\subset\R^n\oplus{\R^n}^*,
\end{equation}
for all $t\in[0,1]$. In formula \eqref{eq:defellt} we think of $g$ as a linear map from $\R^n$ to
${\R^n}^*$; this kind of identification will be made implicitly when necessary in the rest of the paper.
We denote by
$t\mapsto\Phi(t)$ the {\em flow\/} of the Morse--Sturm system
\eqref{eq:MS}, i.e., for every $t\in[0,1]$, $\Phi(t)$ is the unique linear isomorphism of
$\R^n\oplus{\R^n}^*$ such that
\[\Phi(t)\big(J(0),gJ'(0)\big)=\big(J(t),gJ'(t)\big),\]
for every solution $J$ of \eqref{eq:MS}. Observe that $\Phi$ is a $C^1$ curve is the general linear
group of $\R^n\oplus{\R^n}^*$ satisfying the matrix differential equation $\Phi'(t)=X(t)\Phi(t)$ with
initial condition $\Phi(a)=\Id$, where $X$ is given by:
\begin{equation}\label{eq:defX}
X(t)=\begin{pmatrix}0&\hskip-7pt\hphantom{{}^{-1}}g^{-1}\\gR(t)&\hskip-7pt0\end{pmatrix}.
\end{equation}
The $g$-symmetry of $R$ implies that $X$ is a curve in $\spl(2n,\R)$ and hence $\Phi$ is actually a
$C^1$ curve in $\Spl(2n,\R)$. Set $L_0=\{0\}\oplus{\R^n}^*$ and consider the smooth map:
\begin{equation}\label{eq:defbeta}
\beta:\Spl(2n,\R)\longrightarrow\Lambda
\end{equation}
defined by $\beta(\Phi)=\Phi(L_0)$. We have:
\begin{equation}\label{eq:ellbetaPhi}
\ell=\beta\circ\Phi;
\end{equation}
in particular $\ell$ is a $C^1$ curve in the Lagrangian Grassmannian $\Lambda$.

By our construction, conjugate points along $\gamma$ correspond 
to the conjugate instants of the Morse--Sturm system \eqref{eq:MS}, i.e.,
instants $t_0\in\left]0,1\right]$ such that there exists a non zero solution
$J$ of \eqref{eq:MS} with
$J(0)=J(t_0)=0$. Observe that an instant
$t_0\in\left]0,1\right]$ is conjugate iff $\ell(t)$ is {\em not\/} transversal to $L_0$, in which case
the multiplicity of $t_0$ coincides with the dimension of $\ell(t)\cap L_0$. For $k=0,1,\ldots,n$ we
set:
\[\Lambda_k(L_0)=\big\{L\in\Lambda:\mathrm{dim}(L\cap
L_0)=k\big\}\quad\text{and}\quad\Lambda_{\ge1}(L_0)=\bigcup_{k=1}^n\Lambda_k(L_0).\]
Each $\Lambda_k(L_0)$ is a connected real-analytic embedded submanifold of $\Lambda$ having codimension
$\frac12k(k+1)$ in $\Lambda$; the set $\Lambda_{\ge1}(L_0)$ is not a submanifold, but it is a compact
algebraic subvariety of $\Lambda$ whose regular part is $\Lambda_1(L_0)$. The conjugate instants of the
Morse--Sturm system are the instants when
$\ell$ crosses
$\Lambda_{\ge1}(L_0)$. The Maslov index of a curve in $\Lambda$ with endpoints in $\Lambda_0(L_0)$ is
defined as an {\em intersection number\/} of the curve with the algebraic variety
$\Lambda_{\ge1}(L_0)$. The intersection theory needed in this context can for instance be formalized by
an algebraic topological approach. Namely, the first singular relative homology group
$H_1(\Lambda,\Lambda_0(L_0))$ with integer coefficients is infinite cyclic and a generator can be
canonically described in terms of the symplectic form $\omega$.  
\begin{defin}\label{thm:defMaslov}
Let $l:[a,b]\to\Lambda$ be a continuous curve with endpoints in $\Lambda_0(L_0)$. The {\em Maslov
index\/} of $l$, denoted by $\iMaslov(l)$, is the integer number corresponding to the homology class
defined by $l$ in $H_1(\Lambda,\Lambda_0(L_0))$.
\end{defin}

The Maslov index of curves in $\Lambda$ is additive by concatenation,
since the same property holds for the relative homology class.

If $\ell$ is the curve defined in \eqref{eq:defellt} then the initial endpoint $\ell(0)=L_0$ is not in
$\Lambda_0(L_0)$; if $t=1$ is conjugate then a similar problem occur, i.e.,
$\ell(1)\not\in\Lambda_0(L_0)$. However, it is known that there are no conjugate instants in a
neighborhood of $t=0$ and hence we can give the following:
\begin{defin}\label{thm:defMaslSturm}
Assume that $\gamma(1)$ is not conjugate. The {\em Maslov index\/}
of the geodesic $\gamma$, denoted
$\iMaslov(\gamma)$, is defined as the Maslov index of the curve $\ell\vert_{[ \varepsilon,1]}$, where
$\varepsilon>0$ is chosen such that there are no conjugate instants in $\left]0,\varepsilon\right]$.
\end{defin}
The Maslov index of a geodesic can be computed as an algebraic count of the conjugate points. In order
to make this statement precise, let us recall a few more facts about the geometry of the
Lagrangian Grassmannian. For $L\in\Lambda$, there exists a {\em natural\/} identification
\[T_L\Lambda\cong\Bils(L)\]
of the tangent space $T_L\Lambda$ with the space $\Bils(L)$ of symmetric
bilinear forms on
$L$. Given a $C^1$ curve $l:[a,b]\to\Lambda$ we say that $l$ has a {\em nondegenerate intersection\/}
with
$\Lambda_{\ge1}(L_0)$ at $t=t_0$ if $l(t_0)\in\Lambda_{\ge1}(L_0)$ and the symmetric bilinear form
$l'(t_0)$ is nondegenerate on the space $l(t_0)\cap L_0$; in case $l(t_0)\in\Lambda_1(L_0)$ then the
intersection is nondegenerate precisely when it is {\em transversal\/} in the standard sense of
differential topology. Nondegenerate intersections with $\Lambda_{\ge1}(L_0)$ are isolated; in case all
intersections of a $C^1$ curve $l$ with $\Lambda_{\ge1}(L_0)$ are nondegenerate, we have the following
differential topological method to compute the Maslov index:
\begin{teo}\label{thm:topodiff}
Let $l:[a,b]\to\Lambda$ be a $C^1$ curve with endpoints in $\Lambda_0(L_0)$ having only nondegenerate
intersections with $\Lambda_{\ge1}(L_0)$. Then $l$ has only a finite number of intersections with
$\Lambda_{\ge1}(L_0)$ and the Maslov index of $l$ is given by:
\[\iMaslov(l)=\sum_{t\in\left]a,b\right[}\sgn\big(l'(t)\vert_{l(t)\cap L_0}\big).\]
\end{teo}
\begin{proof}
See \cite[Section~3]{MPT}.
\end{proof}
We now want to apply Theorem~\ref{thm:topodiff} to the curve $\ell$ defined in \eqref{eq:defellt}; to
this aim, we first have to compute the derivative of $\ell$. Using local coordinates in $\Lambda$ one
can compute the differential of the map $\beta$ as:
\begin{equation}\label{eq:difbeta}
\dd\beta(\Phi)\cdot
A=\omega(A\Phi^{-1}\cdot,\cdot)\vert_{\Phi(L_0)}\in\Bils(\Phi(L_0)),
\end{equation}
for all $\Phi\in\Spl(2n,\R)$ and
all $A\in T_\Phi\Spl(2n,\R)$.
\begin{teo}\label{thm:calcMaslMS}
If $\gamma(t_0)$ is a nondegenerate (hence isolated) conjugate point
along $\gamma$, $t_0\in\left]0,1\right[$, then for $\varepsilon>0$ small enough:
\[\iMaslov(\gamma\vert_{[0,t_0+\varepsilon]})=\iMaslov(\gamma\vert_{[0,t_0-\varepsilon]})+\sgn(t_0).\]
If $\gamma(1)$ is not conjugate, and if all the conjugate
points along $\gamma$ are nondegenerate, then the Maslov index of
$\gamma$ is given by:
\[\iMaslov(\gamma)=\sum_{t\in\left]0,1\right[}\sgn(t_0).\]
\end{teo}
\begin{proof}
Using the additivity by concatenation of the Maslov index of curves in $\Lambda$,
the result  is an easy
consequence of Theorem~\ref{thm:topodiff},  where formulas
\eqref{eq:defX},
\eqref{eq:ellbetaPhi} and
\eqref{eq:difbeta} are used to compute $\ell'(t)\vert_{\ell(t)\cap L_0}$.
\end{proof}

\subsection{The Maslov index as a relative index}
\label{sub:maslovrel}
We will now relate the Maslov index of a geodesic with the spectral
flow of the path of Fredholm operators obtained from the index form.

Given a geodesic $\gamma:[0,1]\to M$, the {\em index form\/} is the bounded
symmetric bilinear form $I$ defined on the space $\Hcal_\gamma$ of all
vector fields of Sobolev class $H^1$ along $\gamma$ and vanishing at the
endpoints given by:
\[I(V,W)=\int_0^1\Big[g\big(\Ddt V,\Ddt W\big)+g\big(R(\dot\gamma,V)\,\dot\gamma,W\big)\Big]\,\dd t.\]
The index form $I$ is a Fredholm form on $\Hcal_\gamma$ which is realized by a 
strongly indefinite self-adjoint Fredholm operator on $\Hcal_\gamma$ when
$g$ is neither positive nor negative definite.

Set $k=n_-(g)$; a {\em maximal negative distribution along $\gamma$\/}
is a smooth selection $\Delta=(\Delta_t)_{t\in[0,1]}$ of $k$-dimensional 
subspaces of $T_{\gamma(t)}M$
such that $g\vert_{\Delta_t}$ is negative definite for all $t$.
Given a maximal negative distribution $\Delta$ along $\gamma$, denote by $\Scal^\Delta$
the closed subspace of $\Hcal_\gamma$ given by:
\begin{equation}\label{eq:Sdelta}
\Scal^\Delta=\Big\{V\in\Hcal_\gamma:V(t)\in\Delta_t,\ \text{for all
$t\in[0,1]$}\Big\}.\end{equation} The $I$-orthogonal space to $\Scal^\Delta$ has been studied in
\cite{topology}, and it can be characterized as the space of vector fields $V$ along $\gamma$ that
are ``Jacobi in the directions of $\Delta$'', i.e., such that $\Ddtt V-R(\dot\gamma,V)\,\dot\gamma$
is $g$-orthogonal to $\Delta$ pointwise (see \cite[Section~5]{topology}). 
\begin{prop}\label{thm:indexformRCP}
The restriction $I\vert_{\Scal^\Delta}$ is RCPNI and the restriction
$I\vert_{(\Scal^\Delta)^{\perp_I}}$ is RCPPI. Moreover, if $\gamma(1)$ is not conjugate,
the index of $I$ relatively to $\Scal^\Delta$ equals the Maslov index of $\gamma$:
\begin{equation}\label{eq:maslovrel}
\ind_{\Scal^\Delta}(I)=\iMaslov(\gamma).
\end{equation}
\end{prop}
\begin{proof}
The first statement in the thesis is proven in \cite[Proposition~5.25]{topology},
the second statement is proven in \cite[Lemma~2.6.6]{london}.
Equality \eqref{eq:maslovrel} follows from Proposition~\ref{thm:indrelindcoind}
and the semi-Riemannian Morse index theorem \cite[Theorem~5.2]{topology}, that
gives us the equality:
\[\iMaslov(\gamma)=n_-\Big(I\big\vert_{(\Scal^\Delta)^{\perp_I}}\Big)
-n_+\Big(I\big\vert_{\Scal^\Delta}\Big).\qedhere\]
\end{proof}
%%%%%%%%
%%%%%%%%
%%%%%%%%
\end{section}

%%%%%%%%%%%%%%%%
%%%%%%%%%%%%%%%%
\begin{section}{The geometrical bifurcation problem}\label{sec:geometrical}

Let $\gamma:[0,1]\to M$ be a geodesic in $(M,g)$, with $p=\gamma(0)$ and $q=\gamma(1)$; let us
consider again a $g$-orthonormal basis $v_1,\ldots,v_n$ of $T_{\gamma(0)}M$ and assume that
the first $k$ vectors $v_1,\ldots,v_k$ generate a $g$-negative space, while
the $v_{k+1},\ldots,v_n$ generate a $g$-positive space.
Let us consider again the parallel transport of the $v_i$'s along $\gamma$,
that will be denoted by $V_1,\ldots,V_n$. Observe that, since parallel transport is
an isometry, then, for all $t\in[0,1]$, the vectors $V_1(t),\ldots,V_k(t)$ generate a $g$-negative
subspace of $T_{\gamma(t)}M$, that will be denoted by $D^-_t$, and $V_{k+1}(t),\ldots,V_n(t)$ generate
a $g$-positive subspace of $T_{\gamma(t)}M$, denoted by $D^+_t$.

We fix a positive number $\varepsilon_0<1$ such that there are no
conjugate points to $p$ along $\gamma$ in the interval $\left]0,\varepsilon_0\right]$.
Finally, let us
define an auxiliary positive definite inner product on each $T_{\gamma(t)}M$, that will be denoted by
$\gR$, by declaring that the basis
$V_1(t),\ldots,V_n(t)$ be orthonormal.  

\subsection{Reduction to a standard bifurcation problem}
\label{sub:redstbif}
For all $s\in[\varepsilon_0,1]$, let $\Omega_s$ denote
the manifold of all curves $x:[0,s]\to M$ of Sobolev  class $H^1$
such that $x(0)=\gamma(0)=p$ and $x(s)=\gamma(s)$.
It is well l known that $\Omega_s$ has the structure of an infinite dimensional
Hilbert manifold, modeled on the Hilbert space $H^1_0([0,s],\R^n)$.
The geodesic action functional $F_s:\Omega_s\to\R$, defined by:
\begin{equation}\label{eq:defFs}
F_s(x)=\frac12\int_0^sg(\dot x,\dot x)\,\mathrm dt,
\end{equation}
is smooth, and its critical points are precisely the geodesics
in $M$ from $p$ to $\gamma(s)$. For each $x\in\Omega_s$, the
tangent space $T_x\Omega_s$ is identified with the
Hilbertable space:
\[T_x\Omega_s=\big\{V\ \text{vector field along $x$ of class $H^1$}:V(0)=0,\ V(s)=0\big\};\]
we choose the following Hilbert space inner product on each $T_x\Omega_s$:
\begin{equation}\label{eq:innprodTx}
\langle V,W\rangle=\int_0^s\gR\big(\Ddt V,\Ddt W\big)\,\mathrm dt,\quad
V,W\in T_x\Omega_s.\end{equation}
\begin{convention}
In what follows, each tangent space $T_\gamma\Omega_s$ will be identified with the Hilbert
space $H^1_0([0,s],\R^n)$ via the parallel frame $V_1,\ldots V_n$:
\begin{equation}\label{eq:isomorph}
H^1_0([0,s],\R^n)\ni (f_1,\ldots,f_n)\cong \sum_{i=1}^nf_iV_i\in
T_\gamma\Omega_s.\end{equation} 
Since the frame $V_1,\ldots V_n$ is parallel, the semi-Riemannian metric $g$ is
carried  by the isomorphism \eqref{eq:isomorph} into a fixed
symmetric bilinear form $g$ on $\R^n$, covariant
differentiation along
$\gamma$ is carried into standard differentiation of curves in
$\R^n$, and the inner product \eqref{eq:innprodTx} becomes the standard 
$H^1_0$-inner product in $H^1_0([0,s],\R^n)$:
\begin{equation}\label{eq:innproductH}
\phantom{%
V,W\in H^1_0([0,s],\R^n).}\langle V,W\rangle=\int_0^s\gR(V',W')\,\mathrm dt,\quad
V,W\in H^1_0([0,s],\R^n).\end{equation}
Similarly, the subspaces $D^-_t$ and $D^+_t$ of $T_{\gamma(t)}M$
are carried to constant subspaces denoted respectively $D^-$ and $D^+$.
Moreover, the curvature tensor $R$ along $\gamma$ is carried by the isomorphism
\eqref{eq:isomorph} into a smooth curve $t\mapsto R(t)$ of $g$-symmetric
endomorphisms of $\R^n$.
\end{convention}

For $\varepsilon_0\le s_1\le s_2\le 1$ and $x\in\Omega_{s_2}$, there is an obvious isometric
embedding $T_x\Omega_{s_1}\to T_x\Omega_{s_2}$ obtained by extension to $0$
in $\left]s_1,s_2\right]$, but for our purposes we will need 
a deeper identification of (suitable open subsets of) all the Hilbert manifolds $\Omega_s$.
Towards this goal, we do the following construction. 
Let $\rho>0$ be a positive number, assume for the moment that
$\rho$ is less than the injectivity radius of $M$ at $\gamma(s)$ for all $s\in[\varepsilon_0,1]$;
a further restriction for the choice of $\rho$ will be given in what follows.
Let $\mathcal W$ be the open ball of radius $\rho$ centered at $0$ in
$H^1_0([0,1],\R^n)\cong T_\gamma\Omega_1$ and, for all $s\in[\varepsilon_0,1]$, let
$\mathcal W_s$ be the neighborhood of $0$ in $H^1_0([0,s],\R^n)\cong T_\gamma\Omega_s$ given by
the image of $\mathcal W$ by the reparameterization map $\Phi_s$
defined  by:
\begin{equation}\label{eq:defRs} H_0^1([0,1],\R^n)\ni
V\longmapsto V(s^{-1}\cdot)\in H_0^1([0,s],\R^n).
\end{equation}
Finally, for all $s\in[\varepsilon_0,1]$, let $\widetilde {\mathcal W}_s$ be the subset of
$\Omega_s$ obtained as the image of $\mathcal W_s$ by the map:
\[V\longmapsto \mathrm{EXP}(V),\]
where \begin{equation}\label{eq:defEXP}\mathrm{EXP}(V)(t)=\exp_{\gamma(t)}V(t).\end{equation}
Since $\exp_{\gamma(t)}$ is a local diffeomorphism between a neighborhood of
$0$ in $T_{\gamma(t)}M$ and a neighborhood of $\gamma(t)$ in $M$, 
it is easily seen that the positive number $\rho$ above can be chosen small
enough so that, for all $s\in[\varepsilon_0,1]$, $\widetilde{\mathcal W}_s$ is an open subset 
of $\Omega_s$ (containing $\gamma$) and $\mathrm{EXP}$  is a diffeomorphism
between $\mathcal W_s$ and $\widetilde{\mathcal W}_s$.

In conclusion, we have a  family of diffeomorphisms $\Psi_s:\mathcal W\to\widetilde{\mathcal
W}_s$:
\[\Psi_s=\mathrm{EXP}\circ \Phi_s,\]
and we can define a family $(f_s)_{s\in[\varepsilon_0,1]}$ of smooth functionals on
$\mathcal W$ by setting:
\[f_s=F_s\circ\Psi_s;\]
observe that $\Psi_s(0)=\gamma\vert_{[0,s]}$ for all $s$.
\begin{prop}\label{thm:riparamfunct}
$(f_s)_s$ is a smooth family of functionals on $\mathcal W$.
For each $s\in[\varepsilon_0,1]$, a point $x\in \mathcal W$ is
a critical point of $f_s$ if and only if $\Psi_s(x)$ is a geodesic
in $M$ from $p$ to $\gamma(s)$ in $\widetilde{\mathcal W}_s$.
In particular, $0$ is a critical point of $f_s$ for all $s$, and
every  geodesic in $M$ from $p$ to
$\gamma(s)$  sufficiently close to $\gamma$ in the $H^1$-topology 
is obtained from a critical point of $f_s$ in $\mathcal W$.
The second variation of $f_s$ at $0$ is given by the bounded
symmetric bilinear form $I_s$ on $H^1_0([0,1],\R^n)$ defined by:
\begin{equation}\label{eq:secvars}
I_s(V,W)=\int_0^1\Big[\frac1s
g\big(V'(t),W'(t)\big)+sg\big(R(st)V(t),W(t)\big)\Big]\,\mathrm dt.
\end{equation}
\end{prop}
\begin{proof}
The smoothness of $s\mapsto f_s$ follows immediately from the smoothness
of the exponential map and of the reparameterization map $s\mapsto \Phi_s$.
Since $\Psi_s$ is a diffeomorphism for all $s$, the critical
points of $f_s$ are precisely the inverse image through $\Psi_s$ of the
critical points of $F_s$, and the second statement of the thesis is clear 
from our construction. As to the second variation of $f_s$ at $0$,  
formula \eqref{eq:secvars} is easily obtained from the classical
second variation formula for the geodesic action functional $F_s$
at the geodesic $\gamma\vert_{[0,s]}$:
\[\mathrm d^2F_s(\gamma)[V,W]=\int_0^s\Big[g\big(V'(t),W'(t)\big)+g\big(R(t)V(t),W(t)\big)
\Big]\,\mathrm d\tau \]
with the change of variable $t=\tau s^{-1}$.
\end{proof}
Proposition~\ref{thm:riparamfunct} gives us the link between the
notion of bifurcation for a smooth family
of functionals and the geodesic  bifurcation problem
discussed in the introduction.
\subsection{Conjugate points and bifurcation}
We will now compute the spectral flow of the smooth curve of strongly indefinite
self-adjoint Fredholm operators on $H^1_0([0,1],\R^n)$ associated to the curve of symmetric
bilinear forms \eqref{eq:secvars}.
\begin{lem}\label{thm:fredstrind}
For all $s\in[\varepsilon_0,1]$, the bilinear form $I_s$ of \eqref{eq:secvars} is realized by
a bounded self-adjoint Fredholm operator $S_s$ on $H^1_0([0,1],\R^n)$.
If $0<n_-(g)<n$, then $S_s$ is strongly indefinite. If $\gamma(1)$ is not conjugate to
$\gamma(0)$ along $\gamma$, then the endpoints of the path \[[\varepsilon_0,1]\ni s\longmapsto S_s\in
\Fcalsa_*\big(H^1_0[0,1],\R^n)\big)\] are invertible.
\end{lem}
\begin{proof}
The bilinear form $I_s$ in \eqref{eq:secvars} is symmetric and bounded
in the $H^1$-topology, hence $S_s$ is self-adjoint and bounded.

The bilinear form $G$ on $H^1_0([0,1],\R^n)$ defined by $(V,W)\mapsto \tfrac1s\int_0^1g(V',W')\,\mathrm
dt$ is realized by an invertible operator, because $g$ is
nondegenerate. The difference $I_s-G$ is realized by a self-adjoint {\em compact\/}
operator on $H^1_0([0,1],\R^n)$, because it is clearly continuous in the
$C^0$-topology, and the inclusion $H^1_0\hookrightarrow C^0 $ is compact.
This proves that $S_s$ is Fredholm.

Fix now $s\in\left]\varepsilon_0,1\right]$, $s_0\in\left]\varepsilon_0,s\right[$ and, assuming
that $0<n_-(g)<n$, choose $v_+$ and $v_-$ in $\R^n$
with $g(v_+,v_+)>0$ and $g(v_-,v_-)<0$. Let $J_+$ (resp., $J_-$) be the unique Jacobi field
along $\gamma$ such that $J_+(s_0)=v_+$  (resp.,
$J_-(s_0)=v_-$). An easy computations shows that, for all $f\in H^1_0([0,s],\R^n)$,
the following equalities hold:
\[I_s(fJ_+,fJ_+)=\int_0^s(f')^2g(J_+,J_+)\,\mathrm dt,\quad 
I_s(fJ_-,fJ_-)=\int_0^s(f')^2g(J_-,J_-).\]
It follows in particular that $I_s$ is positive definite on the infinite
dimensional subspace of $H^1_0([0,s],\R^n)$ consisting of vector
fields of the form $fJ_+$, with $f$ having a fixed small support 
around $s_0$, and $I_s$ is negative definite on the space
of vector fields of the form $fJ_-$. Hence, $S_s$ is strongly indefinite.

Since $S_s$ is Fredholm of index zero, then $S_s$ is invertible if and only 
it is injective, i.e., if and only if $I_s$ has trivial kernel, that is, if and
only if $\gamma(s)$ is not conjugate to $\gamma(0)$ along $\gamma$. Hence,
the last statement in the thesis comes from the fact that both
$\gamma(\varepsilon_0)$ and $\gamma(1)$ are not conjugate to
$\gamma(0)$ along $\gamma$.
\end{proof}

\begin{lem}\label{thm:continuousext}
The smooth path $\hat I$ of bounded symmetric bilinear forms
$\left]0,1\right]\ni s\mapsto \hat I_s:=s\cdot I_s$ has a continuous
extension to $0$ which is obtained by setting:
\[\hat I_0(V,W)=\int_0^1g(V',W')\,\mathrm dt.\]
For all $s\in[0,1]$, let $\hat S_s$ be the realization of $\hat I_s$ and assume
that $\gamma(1)$ is not conjugate to $\gamma(0)$ along $\gamma$.
The spectral flow of the path $\hat I:[0,1]\to \Fcalsa_*([0,1],\R^n)$
is equal to the spectral flow of the path $S:[\varepsilon_0,1]\to \Fcalsa_*([0,1],\R^n)$.
\end{lem}
\begin{proof}
From \eqref{eq:secvars} we get:
\begin{equation}\label{eq:percalcKs}\hat I_s(V,W)=\int_0^1\Big[ 
g\big(V'(t),W'(t)\big)+s^2g\big(R(st)V(t),W(t)\big)\Big]\,\mathrm dt\end{equation}
for all $s\in\left]0,1\right]$, and   this formula proves immediately
the first statement in the thesis.

The cogredience invariance of $\spfl$ implies that
multiplication by a positive map does not change the spectral flow;
in particular, the spectral flow of $\hat S$ and of $S$ on the interval
$[\varepsilon_0,1]$ coincide.
Since $\hat S_s$ is invertible for all $s\in[0,\varepsilon_0]$, the spectral flow of  $S$ on
$[\varepsilon_0,1]$ coincide with the spectral flow of $\hat S$ on $[0,1]$.
\end{proof}
We are now ready to compute the spectral flow of the path $S$:
\begin{prop}
\label{thm:calcolaflusso}
Assume that $\gamma(1)$ is not conjugate to $\gamma(0)$ along $\gamma$.
Then the spectral flow of the path $S$ is equal to
$-\iMaslov(\gamma)$.
\end{prop}
\begin{proof}
We will compute the spectral flow of the path $\hat S$ on the interval
$[0,1]$; to this aim, we will use part \eqref{itm:2calcspfl} of 
Proposition~\ref{thm:cogrinvformsf}. We will show that $\hat S_s$
has the form $\mathfrak J+K_s$ for all $s\in[0,1]$, where
$\mathfrak J$ is a fixed symmetry of $H^1_0([0,1],\R^n)$ and $K_s$
is a self-adjoint compact operator. Consider the following closed subspaces
of $H^1_0([0,1],\R^n)$:
\[\begin{split}&H^-=\big\{v\in H^1_0([0,1],\R^n):v(t)\in D^-\ \text{for all $t\in[0,1]$}\big\},\\
&H^+=\big\{v\in H^1_0([0,1],\R^n):v(t)\in D^+\ \text{for all $t\in[0,1]$}\big\}.
\end{split}\]
In the language of subsection~\ref{sub:maslovrel}, $D^-$ corresponds to
a maximal negative distribution, and the space $H^-$ corresponds to the
space $\Scal^\Delta$ of \eqref{eq:Sdelta}.

Clearly, $H^1_0([0,1],\R^n)=H^-\oplus H^+$; moreover,
since $D^-$ and $D^+$ are $\gR$-orthogonal, it follows that 
$H^-$ and $H^+$ are orthogonal subspaces with respect to the inner product
\eqref{eq:innproductH}. Set $\mathfrak J=P_+-P_-$, where 
$P_+$ and $P_-$ are the orthogonal projections onto $H^+$ and $H^-$ respectively.
Recalling that $D^-$ and $D^+$ are $g$-orthogonal, and that $g=\gR$ on $D^+$ and $g=-\gR$ on $D^-$,
 we have:
\[\langle \mathfrak JV,W\rangle=\int_0^1g(V',W')\,\mathrm dt, \]
for all $V,W\in H^1_0([0,1],\R^n)$, and thus:
\[\mathfrak J=\hat S_0.\] As we have observed in the proof
of Lemma~\ref{eq:secvars}, the difference $K_s=\hat S_s-\mathfrak J$ is a compact
operator, and it is computed explicitly from \eqref{eq:percalcKs} as:
\[\langle K_sV,W\rangle=s^2\int_0^1g\big(R(st)V(t),W(t)\big)\,\mathrm dt,\quad
V,W\in H^1_0([0,1],\R^n).\]
Clearly, $V^{\scriptscriptstyle-}(\hat S_0)=H^-$.
We can then use formula~\eqref{thm:cogrinvformsf2}, obtaining that the spectral
flow of the path $\hat S$ is given by the relative index:
\[\spfl(\hat S)=-\ind_{H^-}\big(\hat I_1\big)=-\ind_{H^-}\big(I_1\big)=
-\ind_{\Scal^\Delta}(I).\]
The conclusion follows from Proposition~\ref{thm:indexformRCP}.
\end{proof}

\begin{cor}\label{thm:conjptbif}
Assume that $\gamma(t_0)$ is a nondegenerate conjugate point
along $\gamma$. If $\sgn(t_0)\ne0$, then $\gamma(t_0)$ is a bifurcation
point along $\gamma$. More generally, if $0<t_0<t_1\le1$ are non conjugate
instants along $\gamma$, if $\iMaslov\big(\gamma\vert_{[0,t_0]}\big)
\ne \iMaslov\big(\gamma\vert_{[0,t_1]}\big)$ then there exists at least one
bifurcation instant $t_*\in\left]t_0,t_1\right[$.
\end{cor}
\begin{proof}
By the very same argument used in the proof of Proposition~\ref{thm:calcolaflusso},
for all nonconjugate instant $s\in\left]\varepsilon_0,1\right]$ along $\gamma$,
the spectral flow of the path $S$ on the interval $[\varepsilon_0,s]$ equals
the Maslov index $\iMaslov(\gamma\vert_{[0,s]})$. If $t_0$ is a nondegenerate
(hence isolated) conjugate instant, using the additivity by concatenation
of $\spfl$, from Theorem~\ref{thm:calcMaslMS},
for all $\varepsilon>0$ small
enough we then have that the spectral flow of $S$ in the
interval $[t_0-\varepsilon,t_0+\varepsilon]$ is given by:
\[\begin{split}\spfl(S,[t_0-\varepsilon,t_0+\varepsilon])=&\;\spfl(S,[\varepsilon_0,t_0+\varepsilon])-
\spfl(S,[\varepsilon_0,t_0-\varepsilon])\\=&\;-\iMaslov(\gamma\vert_{[0,t_0+\varepsilon]})+
\iMaslov(\gamma\vert_{[0,t_0-\varepsilon]})=-\sgn(t_0).
\end{split}\]
The conclusion follows from Proposition~\ref{thm:FPR} and Proposition~\ref{thm:riparamfunct}.
The proof of the second statement in the thesis is analogous.
\end{proof}

\begin{cor}\label{thm:corRiemcausLor}
If $(M,g)$ is Riemannian, or if $(M,g)$ is Lorentzian and $\gamma$ is causal
(i.e., timelike or lightlike), then every conjugate point along $\gamma$ is a bifurcation
point.
\end{cor}
\begin{proof}
The signature of every conjugate point along a Riemannian manifold
coincides with its multiplicity; the same is true for causal
Lorentzian geodesic. To see this, assume that $\gamma$
is a causal Lorentzian geodesic and $t_0\in\left]0,1\right]$ is a conjugate instant
along $\gamma$; the field $t\dot\gamma(t)$
is in $\mathbb J$, hence $\mathbb J[t_0]^\perp$ is contained
in $\dot\gamma(t_0)^\perp$. If $\gamma$ is timelike, then
$\dot\gamma(t_0)^\perp$ is spacelike, hence $\sgn\big(g\vert_{\mathbb J[t_0]^\perp}\big)=
\Dim\big(\mathbb J[t_0]^\perp\big)=\mul(t_0)$.
If $\gamma$ is lightlike, then $g$ is positive semi-definite on
$\dot\gamma(t_0)^\perp$; to prove that it is positive definite
on $\mathbb J[t_0]^\perp$ it suffices to show that $\dot\gamma(t_0)$
does not belong to $\mathbb J[t_0]^\perp$. To see this, choose a Jacobi
field $J\in\mathbb J$ along $\gamma$ with the property that $\Ddt J(0)$ is
{\em not\/} orthogonal to $\dot\gamma(0)$. It is easily see that
the functions $t\mapsto g\big(J(t),\dot\gamma(t)\big)$ is affine, and
it is zero at $t=0$. If it were $0$ at $t_0$ then it would identically vanish,
which is impossible because its derivative $g\big(\Ddt J(t),\dot\gamma(t)\big)$
does not vanish at $t=0$. It follows that $\dot\gamma(t_0)$ is not orthogonal
to $J(t_0)$, hence $\dot\gamma(t_0)\not\in\mathbb J[t_0]^\perp$.
\end{proof}
\end{section}

%%%%%%%%%%%%%%%%
%%%%%%%%%%%%%%%%
\begin{section}{Final remarks}\label{sec:final}

\subsection{Focal points}
Assume that $\gamma:[0,1]\to M$ is a geodesic in the semi-Riemannian
manifold $(M,g)$, and let
$P\subset M$ be a smooth submanifold with $\gamma(0)\in P$ and $\dot\gamma(0)\in T_{\gamma(0)}P^\perp$.
We will assume that $P$ is nondegenerate at $\gamma(0)$, i.e., that
$g\vert_{T_{\gamma(0)}P}$ is nondegenerate.
Recall that the {\em second fundamental form\/} of $P$ at $\gamma(0)$
in the normal direction $\dot\gamma(0)$ is the symmetric
bilinear form $S^P_{\dot\gamma(0)}:T_{\gamma(0)}P\times T_{\gamma(0)}P\to\R$
given by:
\[S^P_{\dot\gamma(0)}(v,w)=g\big(\nabla_{v}W,\dot\gamma(0)\big),\]
where $W$ is any local extension of $w$ to a vector field in $P$.
A {\em $P$-Jacobi field\/} along $\gamma$ is a Jacobi field $J$
satisfying the initial conditions:
\begin{equation}\label{eq:Pjacobi}
J(0)\in T_{\gamma(0)}P,\quad g\big(\Ddt J(0),\cdot\big)+S^P_{\dot\gamma(0)}\big(J(0),\cdot\big)=0\ 
\text{on $T_{\gamma(0)}P$}.
\end{equation} 
$P$-Jacobi fields are interpreted geometrically as variational vector
fields along $\gamma$ corresponding to variations of $\gamma$ by geodesics
that start orthogonally at $P$. A $P$-focal point
along $\gamma$ is a point $\gamma(t_0)$ for which there exists
a non zero $P$-Jacobi field $J$ such that $J(t_0)=0$.
Observe that the notion of conjugate point coincides with that
of $P$-focal point in the case that $P$ reduces to a single point
of $M$. Theorems~\ref{thm:topodiff} and \ref{thm:calcMaslMS}
hold also in this case, {\em mutatis mutandis}.

The notions of multiplicity and signature of a $P$-focal point, as well as the notion
of nondegeneracy,
are given in perfect analogy with the same notions for conjugate
points (Subsection~\ref{sub:conjugatepts}) by replacing the space
$\mathbb J$ of \eqref{eq:defJbb} with the 
space $\mathbb  J_P$:
\[\mathbb J_P=\big\{J\ \text{solution of \eqref{eq:Jacobieq} satisfying \eqref{eq:Pjacobi}}\big\}.\] 
Also the definition of Maslov index of
$\gamma$ relatively to the initial submanifold $P$,
that will be denoted by $\iMaslov^P(\gamma)$, is analogous
to the definition of Maslov index of a geodesic in the fixed endpoints case
(Subsection~\ref{sub:maslovgeom}). Namely, for the
correct definition Maslov index relative to the initial
submanifold $P$ it suffices to redefine the curve
$\ell$ given in \eqref{eq:defellt} as:
\[\ell(t)=\Big\{\big(J(t),gJ'(t)\big): J\in\mathbb J_P\Big\}\]
and repeat {\em verbatim\/} the definitions in Subsection~\ref{sub:maslovgeom}. 

\begin{defin}\label{thm:defbifP}
A point $\gamma(t_0)$, $t_0\in\left]0,1\right[$, along a geodesic $\gamma:[0,1]\to M$ starting
orthogonally at $P$ is said to be a {\em bifurcation point relatively to the initial
submanifold $P$} (see Figure~\ref{fig:geobif2}) if there exists a sequence $(p_n)_n$ in $P$
converging to $\gamma(0)$, a sequence of normal vectors $N_n\in T_{p_n}P^\perp$
converging to $\dot\gamma(0)$ in the normal bundle $TP^\perp$ (so that 
the geodesic $t\mapsto\exp_{p_n}(tN_n)$
converges to $\gamma$) and a sequence $(t_n)_n$ in $[0,1]$ converging to $t_0$
such that  $\exp_{p_n}(t_n\cdot N_n)$ belongs to $\gamma\big([0,1]\big)$.
\end{defin}
\begin{figure}
\begin{center}
\psfull
\epsfig{file=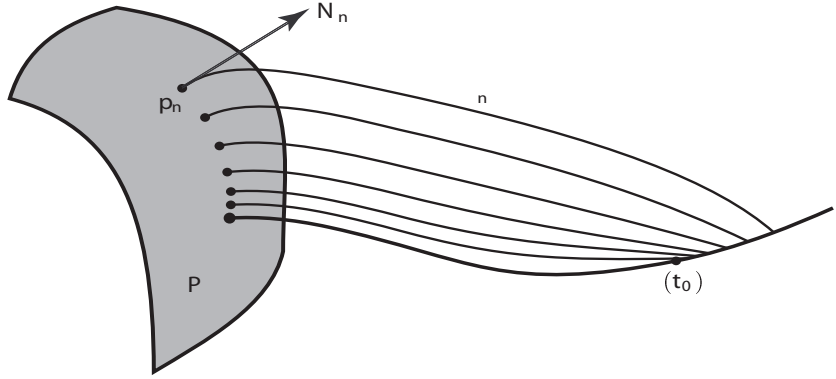,width=11truecm, height=5truecm}
\caption{Bifurcation of geodesics starting orthogonally at a submanifold $P$,
occurring at a $P$-focal point
along $\gamma$.}\label{fig:geobif2}
\end{center}
\end{figure}
The geodesic starting orthogonally at $P$ and terminating at the point
$\gamma(s)$ are critical points of the geodesic action functional
$F_s$ in \eqref{eq:defFs} in the manifold $\Omega^P_s$ of all curves
$x:[0,s]\to M$ of Sobolev class $H^1$ with $x(0)\in P$ and
$x(s)=\gamma(s)$. 
For $x\in \Omega^P_s$, the tangent space $T_x\Omega_s^P$ is 
identified with the space of vector fields $V$ of class $H^1$ along
$x$ such that $V(0)\in T_{x(0)}P$ and $V(s)=0$.
For each $s\in\left]0,1\right]$, the second
variation of $F_s$ at $\gamma\vert_{[0,s]}$ is given
by the symmetric bounded bilinear form $I^P_s$ on $T_\gamma\Omega_s^P$
given by:
\begin{equation}\label{eq:IPs}
I^P_s(V,W)=\int_0^s\Big[g(\Ddt V,\Ddt
W)+g\big(R(\dot\gamma,V)\,\dot\gamma,W\big)
\Big]\,\mathrm d\tau-S^P_{\dot\gamma(0)}\big(V(0),W(0)\big).
\end{equation} 
Using   a parallely transported orthonormal basis
along $\gamma$,
we will identify\footnote{%
Such identification is done in perfect analogy with what discussed in the
Convention on page~\pageref{eq:innprodTx}.}  the tangent space
$T_\gamma\Omega_s^P$ with the Hilbert space 
$H^1_{\mathfrak P}([0,s],\R^n)$ of
all maps 
$V:[0,s]\to \R^n$ of class $H^1$ such that $V(0)\in \mathfrak P$ and 
$V(s)=0$, where $\mathfrak P$ a subspace of $\R^n$ corresponding to
$T_{\gamma(0)}P$ by the above identification of $T_{\gamma(0)}M$ with
$\R^n$, and $\mathcal S$ is the bilinear form on $\mathfrak P$
corresponding to the second fundamental form $S^P_{\dot\gamma(0)}$.
The space $H^1_{\mathfrak P}([0,1],\R^n)$ will be endowed with the
following Hilbert space inner product:
\[\langle V,W\rangle_{\mathfrak P}=\int_0^s\gR(V',W')\,\mathrm dt+\gR\big(V(0),W(0)\big).\]
In order to reduce the focal bifurcation problem to a standard
bifurcation setup, we need to modify slightly the
construction done in Subsection~\ref{sub:redstbif}; this is due to the
fact that the map $\mathrm{EXP}$ as defined in \eqref{eq:defEXP}, when
evaluated on vector fields $V\in T_\gamma\Omega^P_s$, does not produce\footnote{%
Observe indeed that $\exp_{\gamma(0)}v\not\in P$ in general for $v\in T_{\gamma(0)}P$.} 
a curve starting on $P$. However, the reader will quickly convince himself that
the exponential map $\exp_{\gamma(t)}$ in the definition of $\mathrm{EXP}$ in \eqref{eq:defEXP} can be
equivalently  replaced by the exponential map $\widetilde\exp_{\gamma(t)}$
of just about {\em any\/} other metric $\tilde g$ on (an open neighborhood of $\gamma$ in)
$M$. Such replacement will not alter any of the results discussed insofar. 
In order to obtain a well defined map $\mathrm{EXP}$ that sends an open
neighborhood of $0$ in $T_\gamma\Omega_s^P$ diffeomorphically onto an open neighborhood
of $\gamma\vert_{[0,s]}$ in $\Omega_s^P$, it will then suffice to use
the exponential map $\widetilde\exp$ of a (Riemannian) metric
$\tilde g$ defined in an open subset $U\subset M$ containing   $\gamma([0,1])$ 
with the property that $P$ is {\em totally geodesic\/}  relatively to $\tilde g$ near $\gamma(0)$.
Such a metric $\tilde g$ is easily found in a neighborhood of $\gamma(0)$ in $M$ using
a submanifold chart for $P$ around $\gamma(0)$, and then extended using a partition of unity.
Once this has been clarified, the reduction of the focal bifurcation
problem to a standard bifurcation setup is done in perfect analogy
with what discussed in Subsection~\ref{sub:redstbif}: for all $s\in\left]0,1\right]$, 
an open neighborhood $\widetilde{\mathcal W}_s$ of $\gamma\vert_{[0,s]}$ in $\Omega^P_s$ is identified
via $\mathrm{EXP}$ and a reparameterization map with a fixed open neighborhood
$\mathcal W$ of $0$ in $H^1_{\mathfrak P}([0,1],\R^n)$. This identification
carries $\gamma\vert_{[0,s]}$ to $0$ for all $s$, and the family $(F_s)$ of geodesic action
functionals  on $\widetilde{\mathcal W}_s$ to a smooth curve
of functionals $f_s$ on $\mathcal W$. 
For all $s\in\left]0,1\right]$, the second variation of $f_s$ at $0$
is identified with a symmetric bilinear form $I^P_s$
on $H^1_{\mathfrak P}([0,1],\R^n)$ given by:
\begin{equation}\label{eq:hatIP}
I^P_s(V,W)=\int_0^1\Big[\frac1sg\big(V'(t),W'(t)\big)+sg\big(R(st)V(t),W(t)\big)\Big]\,\mathrm dt - 
\mathcal S\big(V(0),W(0)\big).\end{equation} 
The smooth family of bilinear form $\hat I^P_s:=s\cdot I^P_s$, given by:
\[\hat I^P_s(V,W)=\int_0^1\Big[g\big(V'(t),W'(t)\big)+s^2g\big(R(st)V(t),W(t)\big)\Big]\,\mathrm
dt - 
s\mathcal S\big(V(0),W(0)\big)\]
has a continuous extension to $s=0$.

Choose a maximal negative distribution $\Delta$ along $\gamma$
and define the space $\Scal^\Delta$ as in \eqref{eq:Sdelta}; the
semi-Riemannian index theorem \cite[Theorem~5.2]{topology}
tells us that in this case, the $P$-Maslov index $\iMaslov^P(\gamma)$ is
given by:
\begin{equation}\label{eq:iMaslovPindcoind}
\iMaslov^P(\gamma)=n_-\big(I^P_1\vert_{(\Scal^\Delta)^{\perp_{I_1}}}\big)
-n_+\big(I^P_1\vert_{\Scal^\Delta}\big)-n_-\big(g\vert_{T_{\gamma(0)}P}\big),
\end{equation}
where $n_-\big(g\vert_{T_{\gamma(0)}P}\big)$ is the index of the restriction
of $g$ to $T_{\gamma(0)}P$. Recall that this restriction
is assumed nondegenerate, and, by continuity, $g$ will be also nondegenerate
when restricted to tangent spaces of $P$ at points near $\gamma(0)$.
In particular, the index $n_-\big(g\vert_{T_qP}\big)$ is constant
for $q$ near $\gamma(0)$ in $P$.

Using Proposition~\ref{thm:cogrinvformsf} (recall formula~\eqref{thm:cogrinvformsf2}), from
\eqref{eq:iMaslovPindcoind} we get that the spectral flow
of the path $\hat S$ of Fredholm operators realizing the bilinear
form $\hat I^P_s$ in $H^1_{\mathfrak P}([0,1],\R^n)$ with respect to the
inner product $\langle\cdot,\cdot\rangle_{\mathfrak P}$ is given by:
\[\spfl(\hat S)=-\iMaslov^P(\gamma)-n_-\big(g\vert_{T_{\gamma(0)}P}\big).\]

The above construction and  arguments analogous to those used in the
proofs of Corollary~\ref{thm:conjptbif} and Corollary~\ref{thm:corRiemcausLor}
give us the following conclusion:
\begin{prop}\label{thm:focalbif}
Let $(M,g)$ be a semi-Riemannian manifold, $P\subset M$ a smooth submanifold
and $\gamma:[0,1]\to M$ starting orthogonally on $P$; 
assume that $P$ is nondegenerate at $\gamma(0)$.
Then, every non degenerate $P$-focal point with non zero signature is
a bifurcation point relatively to the initial
submanifold $P$. More generally, if $[a,b]\subset \left]0,1\right]$
is such that $\iMaslov^P\big(\gamma\vert_{[0,a]}\big)\ne
\iMaslov^P\big(\gamma\vert_{[0,b]}\big)$, then there exists
at least one bifurcation point relatively to the initial
submanifold $P$ along $\gamma\vert_{\left] a,b\right[}$.

If $(M,g)$ is Riemannian, or if $(M,g)$ is Lorentzian and $\gamma$ is
causal, then every $P$-focal point along $\gamma$ is a bifurcation
point relatively to $P$.\qed
\end{prop}

\subsection{Branching points along geodesics}
A stronger property than bifurcation can be defined for a
point $\gamma(t_0)$ along a semi-Riemannian geodesic $\gamma$
by requiring the existence of a whole homotopy of geodesics 
$\gamma_s$, $s\in I$ where $I\subset\R$ is a right or a left neighborhood
of $t_0$, such that $\gamma_s(a)=\gamma(a)$,  $\gamma_s(s)=\gamma(s)$,
$\gamma_s\ne\gamma$ and $\gamma_s\to\gamma$ as $s\to t_0$.
This is for instance the case of the conjugate point along
a meridian of the paraboloid mentioned in the Introduction.
A point for which such stronger bifurcation property holds
is called a {\em branching point\/} along $\gamma$.
Using a classical Lyapunov-Schmidt reduction and the implicit function
theorem, it is easy to prove that simple (i.e., multiplicity $1$)
nondegenerate conjugate points along geodesics are
branching points.

\subsection{Bifurcation by geodesics with a fixed causal character}
A different bifurcation problem in the context of semi-Riemannian
geodesics may be formulated by requiring that the non trivial branch
of geodesics have a fixed causal character. This is particularly
interesting in the case of lightlike geodesics in Lorentzian
manifolds, where light bifurcation may be used to model the so-called
{\em gravitational lensing\/} problem in General Relativity. We observe here that the
result of Corollary~\ref{thm:corRiemcausLor} does not apply to this 
situation.

\subsection{Bifurcation at an isolated degenerate conjugate point}
As we have observed (\cite{MPT, fechado}), degenerate conjugate points
along a semi-Riemannian geodesic may accumulate; however, when
the metric is real-analytic, an easy argument shows that
conjugate points must necessarily be isolated. 
In the real-analytic case, the result of
Corollary~\ref{thm:conjptbif} can be generalized to
the case of arbitrary conjugate points in  terms of root functions
and partial multiplicities, in
the spirit of \cite{rabier}.

\end{section}

%%%%%%%%%%%%%%%%
%%%%%%%%%%%%%%%%

%%%%%%%%%%%%%%%%%%%%%%%%%%%%%%%%%%%%%%%%%%%%%%%%%%%%%%%%%%%%%%%%%%%%%%%%%%
%%%%%%%%%%%%%%%%%%%%%%%%%%%%%%%%%%%%%%%%%%%%%%%%%%%%%%%%%%%%%%%%%%%%%%%%%%

\end{document}